\documentclass[11pt]{article}
\usepackage{amsmath,amssymb,latexsym,epsfig,subfigure,multirow,graphicx,longtable,rotating}

\newcommand{\q}{\quad}

\newcommand{\na}{\nabla}

\newcommand{\txf}{\textbf}

\newcommand{\de}{\delta}

\newcommand{\veps}{\varepsilon}

\newcommand{\ka}{\kappa}

\newcommand{\ph}{\phi}

\newcommand{\Ga}{\Gamma}

\newcommand{\Om}{\Omega}

\newcommand{\mb}{\mathbb}

\newcommand{\iy}{\infty}
\newcommand{\fr}{\frac}
\newcommand{\pa}{\partial}

\begin{document}

\title{A GPU-accelerated Direct-sum Boundary Integral \\
Poisson-Boltzmann Solver}

\author{Weihua~Geng$^1$\footnote{
Corresponding author. Tel: 1-205-3485302, Fax: 1-205-3487067,
Email: wgeng@as.ua.edu}  ~and Ferosh~Jacob$^2$\\
\\
\small \it     $^1$Department of Mathematics,
University of Alabama, Tuscaloosa, AL 35487, USA\\
\small \it     $^2$Department of Computer Science,
University of Alabama, Tuscaloosa, AL 35487, USA\\
}

\date{\today}
\maketitle

\begin{abstract}
In this paper, we present a GPU-accelerated direct-sum boundary integral method
to solve the linear Poisson-Boltzmann (PB)  equation.
In our method, a well-posed boundary integral formulation
is used to ensure the fast convergence of Krylov subspace based linear algebraic solver such as the GMRES.
The molecular surfaces are discretized with flat triangles and centroid collocation.
To speed up our method, we take advantage of the parallel nature of the boundary integral formulation and parallelize the schemes within CUDA shared memory architecture on GPU.
The schemes use only $11N+6N_c$ size-of-double device memory
 for a biomolecule with $N$ triangular surface elements and $N_c$ partial charges.
Numerical tests of these schemes show well-maintained  accuracy and fast convergence.
The GPU implementation using one GPU card (Nvidia Tesla M2070) achieves 120-150X speed-up to the implementation using one CPU (Intel L5640 2.27GHz). With our approach, solving PB equations on well-discretized molecular surfaces with up to 300,000 boundary elements will take less than about 10 minutes, hence our approach is particularly suitable for fast electrostatics computations on small to medium biomolecules.

\vspace*{1cm}


\end{abstract}

\newpage

\section{Introduction}
Molecular mechanics uses Newton's classical mechanics
to model molecular systems.
The potential energy of all systems in molecular mechanics is calculated using force fields.
Among all components of the force fields,
electrostatics are critical due to their ubiquitous existence
and are expensive to compute since
they are long-range pairwise interactions.
Poisson-Boltzmann (PB) model is an effective alternative
for resolving electrostatics that includes energy, potential and forces
of solvated biomolecules \cite{BakerCurrOp}.
As an implicit solvent approach,
the PB model uses a mean field approximation to trace the solvent effects
and applies Boltzmann distribution to model the mobile ions.
These implicit treatments make the PB model
computationally more efficient compared to explicit solvent models,
in which atomic details of solvent molecules and mobile ions are explicitly described.

In the PB model, the computational domain $\mb{R}^3$
is divided into the solute domain $\Om_1$ and
the solvent domain $\Om_2$ by a closed molecular surface $\Gamma$
such that $\mb{R}^3=\Om_1\cup\Om_2\cup\Gamma$.
The molecular surface $\Gamma$ is formed
by the traces of a spherical solvent probe
rolling in contact with the van del Walls balls of the solute atoms~\cite{Richards, Connolly85}.
The molecule, which is located  in domain $\Om_1$
with dielectric constant $\veps_1$,
is represented
by a set of $N_c$ point charges carrying
$Q_i$ charge in the units of $e_c$, the elementary charge, at positions $\txf{x}_i, i=1, ..., N_c$.
The exterior domain
contains the solvent with dielectric constant $\veps_2$,
as well as mobile ions.
For $\txf{x}=(x, y, z)$,
the PB equation for the electrostatic potential in each domain
is derived from Gauss's law and the Boltzmann distribution. Assuming
weak ionic strength (e.g., the concentration of the physiological saline in a room temperature), the linearized PB equation and its interface jump conditions and boundary conditions have the forms
\begin{align}
\label{eq1}
&\na\cdot(\veps_1(\txf{x})\na\ph_1(\txf{x}))=-\sum_{i=1}^{N_c}
q_i\de (\txf{x}-\txf{x}_i)\q\text{in} ~\,\Om_1,\\
\label{eq2}
&\na\cdot(\veps_2(\txf{x})\na\ph_2(\txf{x}))-\ka^2\ph_2(\txf{x})=0\q
\text{in} ~\,\Om_2,\\
\label{jump}
&\ph_1(\txf{x})=\ph_2(\txf{x}),\q \veps_1
\fr{\pa\ph_1(\txf{x})}{\pa\nu}=\veps_2\fr{\pa\ph_2(\txf{x})}{\pa\nu}\q
\text{on} ~\,\Ga = \partial\Om_1 = \partial\Om_2,\\
\label{radiation}
&\lim_{|\txf{x}|\to\iy}\ph_2(\txf{x})=0,
\end{align}
where $\ph_1$ and $\ph_2$ are the electrostatic potentials in each domain,
$q_i=e_c Q_i/k_B T, i=1, ..., N_c,$
$e_c$ the electron charge,
$k_B$ the Boltzmann's constant,
$T$ the absolute temperature,
$\delta$ the Dirac delta function,
$\ka$ the Debye--H\"{u}ckel parameter, and
$\nu$ the unit outward normal on the interface $\Ga$.
Note the $\kappa^2\phi({\bf x})$ term in Eq.~(\ref{eq2}) is the linearized form of $\kappa^2\sinh(\phi({\bf x}))$ when weak ionic strength is assumed. 

The PB equation is an elliptic equation defined on multiple domains
with discontinuous coefficients across the domain interface.
The PB equation has analytical solutions only for the simple geometries
such as spheres \cite{Kirkwood34} or rods \cite{Holst}.
For molecules with complex geometries, the PB equation can only be solved numerically,
which is challenging due to the non-smoothness of the solution subject to Eq.~(\ref{jump}),
the complex geometry of the interface $\Gamma$,
the singular partial charges in Eq.~(\ref{eq1}), and the boundary conditions
at the infinity in Eq.~(\ref{radiation}).
Many numerical PB solvers were developed  and they can be roughly divided into two categories: 
1) the 3-D mesh-based finite difference/finite element methods
\cite{Rocchia01,Im98,Luo02,bakerSept,Qiao,Chenjcc11};
and 
2) the boundary integral methods
\cite{JBKPB, BFZ, Lu06, yokota-bardhan-knepley-barba-hamada-11,altman-bardhan-white-tidor-09,GengKrasnyjcp12, Bordner, LS, VorSch, bajaj-chen-rand-11}.
All these methods have their own advantages and disadvantages.
For example, the PB solvers embedded in molecular modeling packages,
such as Dephi \cite{Rocchia01}, CHARMM \cite{Im98}, AMBER \cite{Luo02}, and APBS \cite{bakerSept}
use standard seven-point finite difference to discretize the PB equation. Although standard finite difference methods arguably result in reduced accuracy due to the smoothened treatment of Eq.~(\ref{jump}),
the efficient, robust and user-friendly features of these PB solvers
brought their popularity to the computational biophysics/biochemistry community.
Compared with 3-D mesh-based methods, the boundary integral methods have advantages since they
impose the singular partial charge in Eq.~(\ref{eq1}) and the far-field boundary condition in Eq.~(\ref{radiation}) exactly,
use appropriate boundary elements to represent surface geometry to desired precision, and
enforce the continuity condition in Eq.~(\ref{jump}) across the interface explicitly.
They are naturally more convenient methods to attack numerical difficulties
in solving the PB equation.

An often claimed advantage of boundary integral methods
compared to finite difference methods is the reduction of the 3-D differential equation to a 2-D surface integral equation.
However, the finite difference methods generate a 7-band sparse matrix,
while the boundary integral methods constitute a fully dense matrix,
which is prohibitively expensive to store and whose matrix-vector product
is computationally costly.
The remedy is to generate the matrix on-the-fly
and compute the matrix-vector product
with the assistance of fast algorithms for $N$-body problems, such as Fast Multipole Methods (FMM) \cite{BFZ,Lu06,yokota-bardhan-knepley-barba-hamada-11,bajaj-chen-rand-11}
and treecode \cite{GengKrasnyjcp12}. In the last few decades,
the advent of the multicore computers and related parallel architectures
such as MPI and Open MP brought the boundary integral methods
to a superior position to 3-D mesh-based methods.
Recently, the appearance of GPU computing further boosts the boundary element methods, particularly for schemes that have simpler algorithms and use less memory.

GPU computing refers to the use of graphics processing units
for scientific and engineering computing applications
rather than traditionally rendering graphs.
A GPU card can be treated as an array of many simplified CPUs
with reduced but concurrent computing power and more limited but faster memory access to CPUs.
GPUs execute many concurrent threads relatively slowly,
rather than a single thread quickly.
This means that GPU computing is more suitable for problems
with high concurrency, straightforward workflow, low memory requirements,
and infrequent communication.
GPUs have broad areas of application, particularly in speeding up molecular modeling. Molecular dynamics packages such as AMBER \cite{AMBER_GPU_GB, AMBER_GPU_PME} and NAMD \cite{NAMD, NAMD_GPU}
have GPU implementations that achieve significant speedup
compared to CPU implementations.

The matrix-vector product computed in boundary integral formulation
is similar to computing $N$-body problems for particle interactions.
Solving $N$-body problems and related applications
in boundary integral methods are popular targets for GPU computing.
For examples, Nyland et. al~\cite{Nyland} computed $N$-body interactions by direct summation,
Burtscher et.~al~\cite{Burtscher} used the Barnes-Hut treecode~\cite{barnes-hut-86}
to compute the dynamics of 5 million point masses on a 1.3GHz GPU with 240 threads,
obtaining a speedup of 74 over a 2.53 GHz CPU, and
Yokota et.~al~\cite{Yokota} demonstrated a FMM accelerated boundary integral PB solver on 512 GPUs,
achieving 34.6 TFlops.

In this paper, we present a parallel boundary integral PB solver on GPUs.
We adopt the well-posed formulation from Juffer et al.~\cite{JBKPB}, rather than the straightforward integral formulation presented in \cite{Yoon} and applied in \cite{altman-bardhan-white-tidor-09, Yokota}. The N-body summation is computed directly, therefore the scheme is implementation convenient and memory saving for GPU computing. 
The rest of the paper is organized as follows.
In section 2, we provide our algorithms including the well-posed boundary integral formulation
and its discretization, followed by CUDA implementation on GPUs.
In section 3, we present the numerical results,
first on the spherical cavities with multiple partial changes,
whose analytical solutions are available  and then
on a series of proteins with various sizes and geometries.
This paper ends with a section of concluding remarks.

\section{Methods}

We use the well-posed boundary integral formulation
from Juffer's work \cite{JBKPB} together with a flat triangulation
and a centroid collocation.
A factor affecting the accuracy of boundary integral method is the discretization of the surface.
We use a non-uniformed triangular surface from MSMS \cite{Sanner}.
Throughout this paper, we call our GPU-accelerated boundary integral Poisson-Boltzmann solver as GABI-PB solver.

\subsection{Well-posed integral formulation}
The differential PB equations ~(\ref{eq1}) and~(\ref{eq2}) can be converted to boundary integral equations. In order to do this, we first define the fundamental solutions to the Poisson equation~(\ref{eq1}) in $\Omega_1$ and the fundamental solution to the PB equation~(\ref{eq2}) in $\Omega_2$ as
\begin{equation}
G_0(\txf{x}, \txf{y}) = \fr{1}{4\pi|\txf{x}-\txf{y}|},
\quad
\label{eq_potential}
G_\kappa(\txf{x}, \txf{y}) = \fr{e^{-\kappa|\txf{x}-\txf{y}|}}{4\pi|\txf{x}-\txf{y}|}.
\end{equation}
Note $G_0(\txf{x},\txf{y})$ and $G_\kappa(\txf{x}, \txf{y})$ are called
Coulomb and screened Coulomb potentials in electrostatic theory.   
By applying Green's second theorem, and canceling the normal derivative terms with interface conditions in Eq.~(\ref{jump}), the coupled integral equations can be derived as \cite{Yoon}:
\begin{align}
\label{eqbim_1}
\ph_1(\txf{x})=&\int_{\Ga} \left[G_0(\txf{x},
\txf{y})\fr{\pa\ph_1(\txf{y})}{\pa\nu_{\txf{y}}} -\fr{\pa
G_0(\txf{x}, \txf{y})}{\pa\nu_{\txf{y}}} \ph_1(\txf{y})
\right]dS_{\txf{y}}+ \sum_{k=1}^{N_{c}}
q_k G_0(\txf{x}, \txf{y}_k) ,\q\hskip 0.1in\txf{x}\in\Om_1,\\
\label{eqbim_2}
\ph_2(\txf{x})=&\int_{\Ga}\left[-G_\kappa(\txf{x},
\txf{y})\fr{\pa\ph_2(\txf{y})}{\pa\nu_{\txf{y}}}+\fr{\pa
G_\kappa(\txf{x}, \txf{y})}{\pa\nu_{\txf{y}}}
\ph_2(\txf{y})\right]dS_{\txf{y}}, \q\hskip 1.18in \txf{x}\in\Om_2.
\end{align}
However,
straightforward discretization of Eqs.~(\ref{eqbim_1}) and (\ref{eqbim_2}) yields an ill-conditioned linear system, whose condition number dramatically increases as number of boundary elements increases~\cite{Lu07}. Juffer et al. derived a well-posed boundary integral formulation by going through the differentiation of the single-layer potentials and the double-layer potentials \cite{JBKPB}. Here the single-layer potentials are from the induced point charge distributions $G_0$ and $G_{\kappa}$ on surface $\Gamma$, while the double-layer potential are from the induced dipole charge distributions, which are the normal derivatives of $G_0$ and $G_{\kappa}$ on surface $\Gamma$. The desired integral forms are as:
\begin{align}
\label{eqbim_3}
\fr{1}{2}\left(1+\veps\right)\ph_1(\txf{x})&=
\int_{\Ga}\left[K_1(\txf{x}, \txf{y})\fr{\pa\ph_1(\txf{y})}
{\pa\nu_{\txf{y}}}+K_2(\txf{x}, \txf{y})\ph_1(\txf{y})\right]dS_{\txf{y}}+S_{1}(\txf{x}),
\qquad \txf{x}\in\Ga, \\
\label{eqbim_4}
\fr{1}{2}\left(1+\fr{1}{\veps}\right)\fr{\pa\ph_1(\txf{x})}{\pa\nu_{\txf{x}}}&=
\int_{\Ga}\left[K_3(\txf{x}, \txf{y})\fr{\pa\ph_1(\txf{y})}
{\pa\nu_{\txf{y}}}+K_4(\txf{x}, \txf{y})\ph_1(\txf{y})\right]dS_{\txf{y}}
+S_{2}(\txf{x}),
\qquad \txf{x}\in\Ga,
\end{align}
with the notation for the kernels
\begin{align}\nonumber
K_1(\txf{x}, \txf{y})=&\,{G_{0}(\txf{x},\txf{y})}-{G_{\kappa}(\txf{x},\txf{y})}, \hskip 0.75in
K_2(\txf{x}, \txf{y})=\veps\fr{\pa G_{\kappa}(\txf{x},\txf{y})}{\pa\nu_{\txf{y}}}-\fr{\pa
G_{0}(\txf{x},\txf{y})}{\pa\nu_{\txf{y}}},\\
\label{Eq_Ls}
K_3(\txf{x}, \txf{y})=&\,\fr{\pa
G_{0}(\txf{x},\txf{y})}{\pa\nu_{\txf{x}}}-\fr{1}{\veps}\fr{\pa
G_{\kappa}(\txf{x},\txf{y})}{\pa\nu_{\txf{x}}}, \hskip 0.4in
K_4(\txf{x}, \txf{y})=\fr{\pa^2
G_{\ka}(\txf{x},\txf{y})}{\pa\nu_{\txf{x}}\pa\nu_{\txf{y}}}-\fr{\pa^2
G_{0}(\txf{x},\txf{y})}{\pa\nu_{\txf{x}}\pa\nu_{\txf{y}}},\\
S_{1}(\txf{x})=&\,\sum_{k=1}^{N_{c}}q_kG_{0}(\txf{x}, \txf{y}_k), \hskip 1.2in
S_{2}(\txf{x})=\sum_{k=1}^{N_{c}}q_k
\fr{\pa
G_{0}(\txf{x},\txf{y}_k)}{\pa\nu_{\txf{x}}}\label{Eq_Ss}.
\end{align}
and $\veps={\veps_{1}}/{\veps_{2}}$.
Note this is the well-posed Fredholm second kind of integral equation,  which is also our choice in this paper.

\subsection{Discretization}
We discretize the integral equations (\ref{eqbim_3}) and (\ref{eqbim_4})
with the flat triangle and the centroid collocation
(the quadrature point is located at the center of each triangle).
This scheme also assumes that the potential and its normal derivative,
as well as the kernel functions are uniform on each triangle.
When the singularities in kernels occur (${\bf x =y}$), the contribution of this triangle in the integral is then simply removed. This scheme, which provides the convenience of incorporating fast algorithms, such as FMM \cite{Lu07} and treecode \cite{GengKrasnyjcp12}, is in fact widely used in the latest boundary integral methods in solving PB equations.  

In this paper, suppose the triangulation program MSMS \cite{Sanner} discretizes the molecular surface to
a set of $N$ triangular elements connected with $N_v$ vertices and $N_e$ edges.
We then have the relation $N_v+N-N_e=2$ called Euler's polyhedron formula. With potential $\phi_1({\bf x})$ and its normal derivative $\fr{\pa\ph_1(\txf{x})}{\pa\nu_{\txf{x}}}$ at the centroid of each triangular element as the unknowns, Eqs.~(\ref{eqbim_3}) and~(\ref{eqbim_4}) are converted to a linear algebraic system ${\bf Au}={\bf b}$, whose elements are specified as follows. Note we do not explicitly express $\bf A$ in an iterative method as we only concern ${\bf Au}$ on each iteration. 

For $i=1,2,...,N$, the $i$th and the $(i+N)$th element
of the discretized matrix-vector product ${\bf Au}$ are given as
\begin{eqnarray}\label{eqdbef1}
\{{\bf Au}\}_{i} &=& \fr{1}{2}\left(1+\veps\right)\ph_1(\txf{x}_{i})-\sum\limits_{j=1,j\ne i}^{N} W_{j}\left[K_1(\txf{x}_{i}, \txf{x}_{j})\fr{\pa\ph_1(\txf{x}_{j})}{\pa\nu_{\txf{x}_{j}}}+K_2(\txf{x}_{i},\txf{x}_{j})\ph_1(\txf{x}_{j})\right] \\
   \{{\bf Au}\}_{i+N}&=&\fr{1}{2}\left(1+\fr{1}{\veps}\right)\fr{\pa\ph_1(\txf{x}_{i})}{\pa\nu_{\txf{x}_{i}}}-\sum\limits_{j=1,j\ne i}^{N} W_{j}\left[K_3(\txf{x}_{i}, \txf{x}_{j})\fr{\pa\ph_1(\txf{x}_{j})}
{\pa\nu_{\txf{x}_{j}}}+K_4(\txf{x}_{i}, \txf{x}_{j})\ph_1(\txf{x}_{j})\right]\label{eqdbef2}
\end{eqnarray}
where $W_j$ is the area of the $j$th element. 
Note we use ${\bf x}_i$ and ${\bf x}_j$  (instead of ${\bf y}_j$) in all kernels to indicate sources and targets are the same set of points. 
The expressions of ${\bf b}_i$ and ${\bf b}_{i+N}$ are directly obtained from $S_1$ and $S_2$ in Eq.~({\ref{Eq_Ss}). In producing the results of this paper, 
we apply GMRES solver  \cite{Saad}  to solve the linear algebraic system from the discretization of Eqs.~(\ref{eqbim_3}) and (\ref{eqbim_4}), which requires computing the matrix-vector product in Eqs.~(\ref{eqdbef1}) and (\ref{eqdbef2}) on each iteration.


\subsection{Electrostatic solvation energy formulation}
To perform the solvation analysis of an interested biomolecules, the 
electrostatic solvation energy is computed by
\begin{equation}
E_{\rm sol} = \frac{1}{2}\sum_{k=1}^{N_c}q_k\phi_{\rm reac}(\txf{x}_k) =
\frac{1}{2}\displaystyle \sum\limits_{k=1}^{N_c} q_k
\int_\Gamma \left[K_1(\txf{x}_k, \txf{y})\fr{\pa\ph_1(\txf{y})}
{\pa\nu_\txf{y}}+K_2(\txf{x}_k, \txf{y})\ph_1(\txf{y})\right]dS_{\txf{y}},
\label{solvation_energy}
\end{equation}
where $\phi_{\rm reac}({\bf x}_k) = \phi_1({\bf x}_k)-S_1({\bf x}_k)$, whose formulation is the integral part of Eq.~(\ref{solvation_energy}),
is the reaction potential at the $k$th solute atom.
The electrostatic solvation energy,  which can be regarded as the atomistic charge $q_k$ weighted average of the reaction potential $\phi_{\rm reac}$, can effectively characterize the accuracy of a PB solver.

\subsection{GPU/CUDA implementation}
\begin{table}[!t]
\caption{
Pseudocode for GABI-PB solver using GPU}
\begin{center}
\vskip -7.5pt
\fbox{
\begin{tabular}{rl}
1 & On host (CPU)\\
2 & \qquad read biomolecule data (charge and structure) \\
3 & \qquad call MSMS to generate triangulation \\
4 & \qquad copy biomolecule data and triangulation to device \\
5 & On device (GPU)\\
6&  \qquad each thread concurrently computes and stores source terms for assigned triangles\\
7 & \qquad copy source terms on device to host\\
8 & On host \\
9 & \qquad set initial guess ${\bf x}_0$ for GMRES iteration and copy it to device \\
10 & On device\\
11 & \qquad each thread concurrently computes assigned segment of matrix-vector product ${\bf y}={\bf Ax}$ \\
12 & \qquad copy the computed matrix-vector ${\bf y}$ to host memory \\
13 & On host\\
14 & \qquad test for GMRES convergence \\
15 & \qquad\qquad if no, generate new ${\bf x}$ and copy it to device, go to step 10 for the next iteration \\
16 & \qquad\qquad if yes, generate and copy the final solution to device and go to step 17 \\
17 & On device\\
18 & \qquad compute assigned segment of electrostatic solvation energy \\
19 & \qquad copy results in step 18 to host \\
20 & On host\\
21 & \qquad add segments of electrostatic solvation energy and output result \\
\end{tabular}
}\label{tb_GPU}
\end{center}
\label{pseudocode}
\end{table}

The GABI-PB solver uses the boundary integral formulation, which can be conveniently parallelized.
The majority of the computing time is taken by the following subroutines.\\
(1) Compute the source term in Eqs.~(\ref{eqbim_3}) and~(\ref{eqbim_4}). \\
(2) Perform matrix-vector product as in Eqs.~(\ref{eqdbef1}) and (\ref{eqdbef2}) on each GMRES iteration.\\
(3) Compute electrostatic solvation energy.\\
Among all of these subroutines, subroutine (2) is the most expensive one and is repeatedly computed on each GMRES iteration. We compute all of these routines in parallel to minimize the computation time.  Table \ref{tb_GPU}  provides the pseudo code for the overall computation.

In this psuedo code, we divide all the operations to those on host performed by the CPU
and those on device performed by the GPU. We use $N_c$ to denote the number of atoms of the biomolecule and $N$ to represent the number of triangular elements. In the following description, we include the memory use in a pair of parentheses following each variable we would have claimed.  

We first read in the atomistic coordinates (3$N_c$), radius ($N_c$) and charges ($N_c$) of the given biomolecule (step 2)
and call the MSMS \cite{Sanner} program to generate triangulation (step 3) including
the faces and the vertices of the triangulation. 
We convert these triangulation information to 
centroid coordinates (3$N$), normal direction vectors (3$N$) and triangle areas ($N$), and copy them together with the atomistic coordinates, radii and charges to the device (step 4).
So far, the device memory use is $5N_c+7N$ size-of-double.
We then use multithreads on device to compute the source terms (2$N$, but deallocated right after being copied to host)
in Eqs.~(\ref{eqbim_3}) and~(\ref{eqbim_4}) (step 6) and copy them to the host (step 7).
Next, we set the initial guess ${\bf x}_0$ of the solution (currently a zero vector)
and copy it  to device (step 9). Note that, when the GMRES is restarted
(usually this needs to be done every 10-20 steps
to ensure the small number of terms in Krylov subspace expansion),
we use the approximated solution achieved
from the previous iteration as the initial guess ${\bf x}_0$.

What follows is the major computing part. We compute the matrix-vector product ${\bf y}={\bf Ax}$ (step 11) on multi-threads GPU device,
and copy the resulting $\bf y$ to the host (step 12).
These two steps take additionally $4N$ size-of-double device memory.
Then on host we test if the obtained $\bf y$ satisfies
the GMRES convergence criterion (step 14) to decide
whether to start the next iteration (step 15)
or terminate the GMRES
and generate the final solution (step 16). After GMRES convergence criterion is
satisfied, we need to compute the atom-wise components
of the electrostatic solvation energy on device using the final solution copied from the host (step 18).
This step takes additional $N_c$ size-of-double device memory. 
Finally, with the components of the electrostatic solvation energy copied from device (step 19),
we compute the total electrostatic solvation energy on host (step 21).
All together, we have allocated $6N_c+11N$ size-of-double device memory.

In implementing algorithms, we tune the CUDA code with some optimization strategies like loop unrolling, coalesced memory access (using double3 data type instead of three double variables for
some structures like position and charge), and fast CUDA operators like ``rsqrt". We also test and choose the optimized value 256 as the number of thread per block. Although these strategies did not substantially improve the performance, they contribute to the overall performance.

\section{Results}
In this section, we present the numerical results.
We first solve the PB equation
on a spherical cavity with multiple charges at different locations.
The analytical solutions in terms of spherical harmonics
are available \cite{Kirkwood34}.
We then solve the PB equation on a series of 24 proteins
with different sizes and geometries.
These protein structures are downloaded from protein data bank (PDB),
and their charges and hydrogens are added with CHARMM22 force field \cite{MacKerell-CHARMM22}.
We report the electrostatic solvation energy results and
its associated execution time for solving the PB equation and computing electrostatic solvation energy on these proteins
with and without GPU acceleration to demonstrate
the improved efficiency from using GPU.
All algorithms are written in C and CUDA
and compiled by gcc with flag ``-O3"
and nvcc with flag ``-O3 -arch=sm\_20 -use\_fast\_math".
The simulations are performed
with a single CPU (Intel(R) Xeon(R) CPU E5640 @ 2.27GHz with 2G Memory) on a 12-core workstation and one GPU (Nvidia Tesla M2070) card. 
Before we reveal the numerical results, we define order and error.

\subsection{Order and Error}
In this paper, we report the relative $L_{\infty}$ error of the surface potential, which is defined as
\begin{eqnarray}\label{eqerr1}
  e_{\phi}&=& \frac{\max\limits_{i=1,...,N}|\phi^{num}(x_{i})-\phi^{exa}(x_{i})|}{\max\limits_{i=1,...,N}|\phi^{exa}(x_{i})|}\label{eqerror}
\end{eqnarray}
where $N$ is the number of unknowns, also the number of triangular elements of a particular discretized molecular surface.
The notation $\phi^{num}$ represents numerically solved surface potential
and $\phi^{exa}$ denotes the analytical solutions obtained
by Kirkwood's spherical harmonic expansion \cite{Kirkwood34}.
The discretization on the molecular surface has a parameter ``density",
number of vertices per \AA$^2$.

The numerical order of accuracy is computed with
\begin{eqnarray}\label{eqeorder}
  \texttt{order}=\displaystyle\log_{\frac{\texttt{coarse\_mesh}}{\texttt{fine\_mesh}}}\frac{e_{\phi}^{\texttt{coarse}}}{e_{\phi}^{\texttt{fine}}}
\end{eqnarray}
following the convention of numerical analysis, where ``mesh" refers to density for boundary integral methods at both coarse and fine levels.

\subsection{Accuracy tests on a spherical cavity}

\begin{figure}[!thbp]
\begin{center}
\includegraphics[width=2.5in]{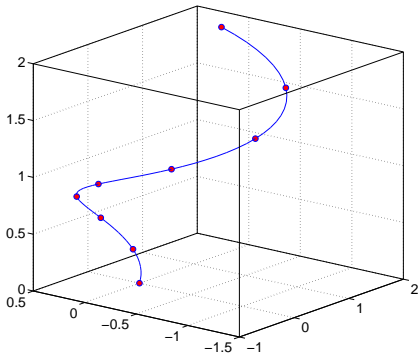}
\caption{Locations of the partial charges inside a spherical cavity with radius $r=4$\AA}
\label{fig_sphere}
\end{center}
\end{figure}

We first perform the numerical tests for the accuracy of the algorithm
on a spherical cavity, where the analytical solution
in terms of spherical harmonics expansion is available \cite{Kirkwood34}.
The test case we designed is a sphere of radius $r=4$\AA
~containing nine partial charges, which are located
along the space curve ${\bf r}(t)=<\frac{3}{4\pi}\cos{t},\frac{3}{4\pi}\sin{t}, \frac{t}{\pi}>$ for $t=0,\frac{\pi}{2},\ldots, 2\pi$. These nine point charges carry $0.1,0.2,\ldots,0.9$ electric charges respectively in the units of $e_c$, the elementary charge. We plot the charge locations in Fig.~\ref{fig_sphere}, which resembles a segment of a biological helix. Table \ref{tb_sphere} reports the numerical results.
In the first column of the table, we increase the MSMS input parameter ``density" (number of vertices per \AA$^2$) by doubling its current value each time. The number of triangular elements are also approximately doubled each time as seen in column 2. The electrostatic solvation energy are reported in column 3 and we can see a consistent pattern that these values are approaching the true value -952.52kcal/mol. Column 4 is
the relative $L_\infty$ error of the surface potential, whose convergence pattern can be better seen from column 5 in terms of orders.
Note that the 0.5th order observed here is relative to the \underbar{area} of the triangular element. If the 1-D \underbar{length} is considered, the order should be about one.

\begin{table}[!htdp]
{\small
\caption{\small Accuracy tests on a spherical cavity: radius $r$=4\AA, nine charges are located on the space curve $<\frac{3}{4\pi}\cos{t},\frac{3}{4\pi}\sin{t}, \frac{t}{\pi}>$ for $t=0:\frac{\pi}{2}:2\pi$; $e_\phi$ is the relative $L_{\infty}$ error in surface potential $\phi$; order is relative to the area of the elements; the exact value of $E_{sol}$ is -952.52 kcal/mol.}
\begin{center}
\begin{tabular}{ccccc}
\hline\hline
density& \# of ele.&$E_\text{sol}$&	$e_\phi$	& order \\\hline
1		&370		&-971.92		&4.20E-02	&-\\	
2		&736		&-968.05		&2.80E-02	&0.58\\
4		&1572		&-964.03		&1.85E-02	&0.60\\
8		&3124		&-961.08		&1.28E-02	&0.53\\
16		&6308		&-958.75		&8.90E-03	&0.52\\
32		&12772	&-956.99		&6.24E-03	&0.51\\
64		&25494	&-955.74		&4.35E-03	&0.52\\
128	&51204	&-954.81		&3.06E-03	&0.51\\
\hline\hline
\end{tabular}
\end{center}
\label{tb_sphere}
}
\end{table}%

\subsection{Accuracy and efficiency test on proteins}
We next solve the PB equation and compute the electrostatic solvation energy on a series of
proteins with different sizes and geometries. The numerical results are reported in Table~\ref{tb_24p}.
In this table, the first column is the index for the convenience of identification. The second row is the four-digit protein data bank (PDB) ID of the corresponding protein. Column 3 is the number of triangular elements when the molecular surfaces of these proteins are discretized. We uniformly choose ``density=10" so that proteins with larger molecular surface areas as seen in column 5 will normally generate larger number of elements. A few exceptions happen when MSMS modifies the given density to fit its triangulation needs, resulting in a slightly mismatched order for the data in columns 3 and 5. Column 4 gives the number of atoms of
associated proteins. We can see in most of the cases the larger number of atoms results in the larger molecular surface area and number of elements.  Column 5 lists the number of iterations. From this column, we noticed the GABI-PB solver converges within 20 steps on the majority of the proteins and the number of iterations does not increase with the increment of the number of elements, which contributes to the well-posed integral formulation. Occasionally, the number of iteration is high, for example for the 7th, 8th, and 24th proteins. This largely attributes
the fact that the MSMS triangulations of these proteins have triangles with very small areas or very short sides. Currently we are
working on a project to improve the triangulation. 
We also plot the number of iterations vs the number of elements in Fig.~\ref{fig_24p}(a) on which the fact that most numbers of iterations are consistently low with a few exceptions is visualized. Column 7 reports the solvation energies. We compare these values with results from our previous work in \cite{Gengjcp07} and consistency of results are observed. Column 8 shows the time for solving PB equation and computing the electrostatic solvation energy on one CPU ($T_1$ in seconds), Column 9 displays the time for that with GPU acceleration ($T_p$ in seconds), and Column 10 reports the ratio between them. From column 9, we can see most of the jobs are finished within 2-3 minutes except jobs on proteins with slow convergence (e.g. the 7th, 8th and 24th proteins). The ratio in the last column demonstrates the overall 120-150X parallel speedup.

\begin{table}[htdp]
{\small
\caption{\small Numerical results for computing the electrostatic solvation energy of 24 proteins: $N$ is the number of elements, $N_c$ is the number of atoms, $N_{it}$ is the number of iterations, area is the molecular surface area, $E_{sol}$ is the electrostatic solvation energy, $T_1$ is the running time on one CPU, and $T_p$ is the running time with GPU acceleration.}
\begin{center}
\begin{tabular}{cccccccccc}
\hline\hline
ID	&PDB	& N	&$N_c$&area(\AA$^2$)&$N_{it}$	&$E_\text{sol}$(kcal/mol)	&CPU $T_1$(s)	&GPU $T_p$(s)& $T_1/T_p$\\
\hline
1	&1ajj	&40496	&519	&2176	&8	&-1145.76	&1323	&11	&125\\
2	&2erl	&43214	&573	&2329	&8	&-961.68	&1410	&11	&124\\
3	&1cbn	&44367	&648	&2377	&8	&-307.03	&1486	&12	&129\\
4	&1vii	&47070	&596	&2488	&11	&-915.18	&2501	&19	&129\\
5	&1fca	&47461	&729	&2558	&8	&-1215.34	&1699	&13	&127\\
6	&1bbl	&49071	&576	&2599	&10	&-1001.47	&2267	&17	&134\\
7	&2pde	&50518	&667	&2727	&99	&-824.91	&25985	&194	&134\\
8	&1sh1	&51186	&702	&2756	&100	&-760.71	&27160	&201	&135\\
9	&1vjw	&52536	&828	&2799	&8	&-1255.93	&2085	&16	&133\\
10	&1uxc	&53602	&809	&2848	&9	&-1154.60	&2437	&18	&136\\
11	&1ptq	&54260	&795	&2910	&10	&-883.77	&2778	&21	&130\\
12	&1bor	&54629	&832	&2910	&12	&-863.91	&3657	&28	&132\\
13	&1fxd	&54692	&824	&2935	&7	&-3344.56	&1963	&15	&129\\
14	&1r69	&57646	&997	&3068	&9	&-1100.83	&2823	&22	&131\\
15	&1mbg	&58473	&903	&3085	&9	&-1368.58	&2901	&22	&132\\
16	&1bpi	&60600	&898	&3247	&24	&-1320.89	&9004	&64	&141\\
17	&1hpt	&61164	&858	&3277	&11	&-825.59	&4248	&32	&133\\
18	&451c	&79202	&1216	&4778	&19	&-1038.20	&11812	&87	&136\\
19	&1svr	&81198	&1435	&4666	&10	&-1731.54	&7294	&53	&138\\
20	&1frd	&81972	&1478	&4387	&9	&-2890.28	&5676	&41	&139\\
21	&1a2s	&84527	&1272	&4457	&16	&-1941.37	&11461	&83	&139\\
22	&1neq	&89457	&1187	&4738	&19	&-1756.02	&15077	&106	&142\\
23	&1a63	&132134	&2065	&7003	&11	&-2404.07	&19804	&139	&143\\
24	&1a7m	&147121	&2809	&7769	&54	&-2184.05	&124326	&852	&146\\
\hline
\end{tabular}
\end{center}
\label{tb_24p}
}
\end{table}

We also plot the results of the last three columns in Fig.~\ref{fig_24p}. From Fig.~\ref{fig_24p}(b),
we see the CPU time $T_1$ is consistently a scale of 100+ of the GPU time $T_p$. Time increases
as the number of elements increases normally. Fig.~\ref{fig_24p}(c) enables us to better observe that the speed-up increases as the number of elements increases, which is advantageous for solving problems of larger sizes.

\begin{figure}[htbp]
\begin{center}
\flushleft \hskip 0.2in (a) \hskip 2.0in (b) \hskip 2.0in (c)\\
\includegraphics[width=2.18in]{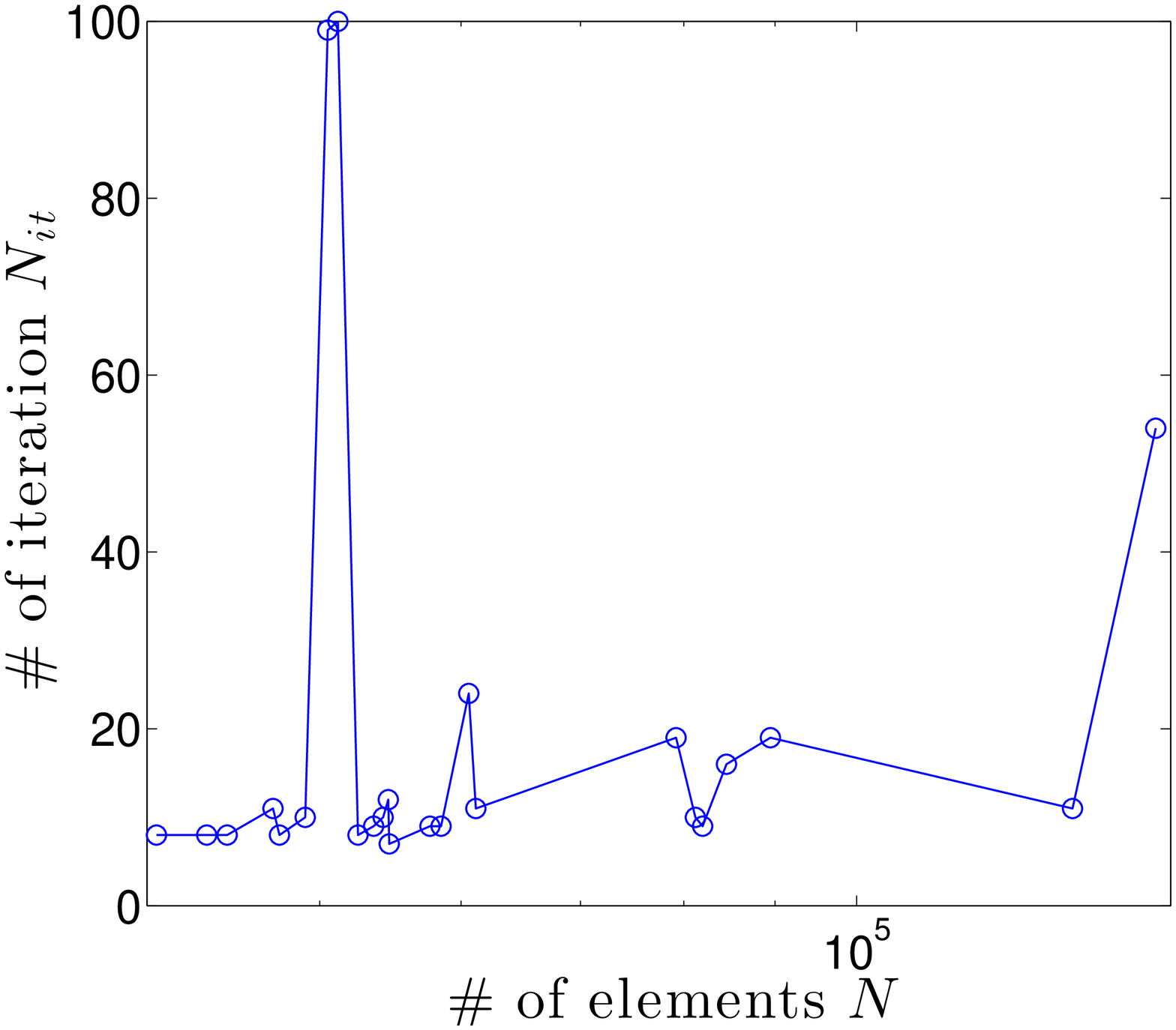}
\includegraphics[width=2.18in]{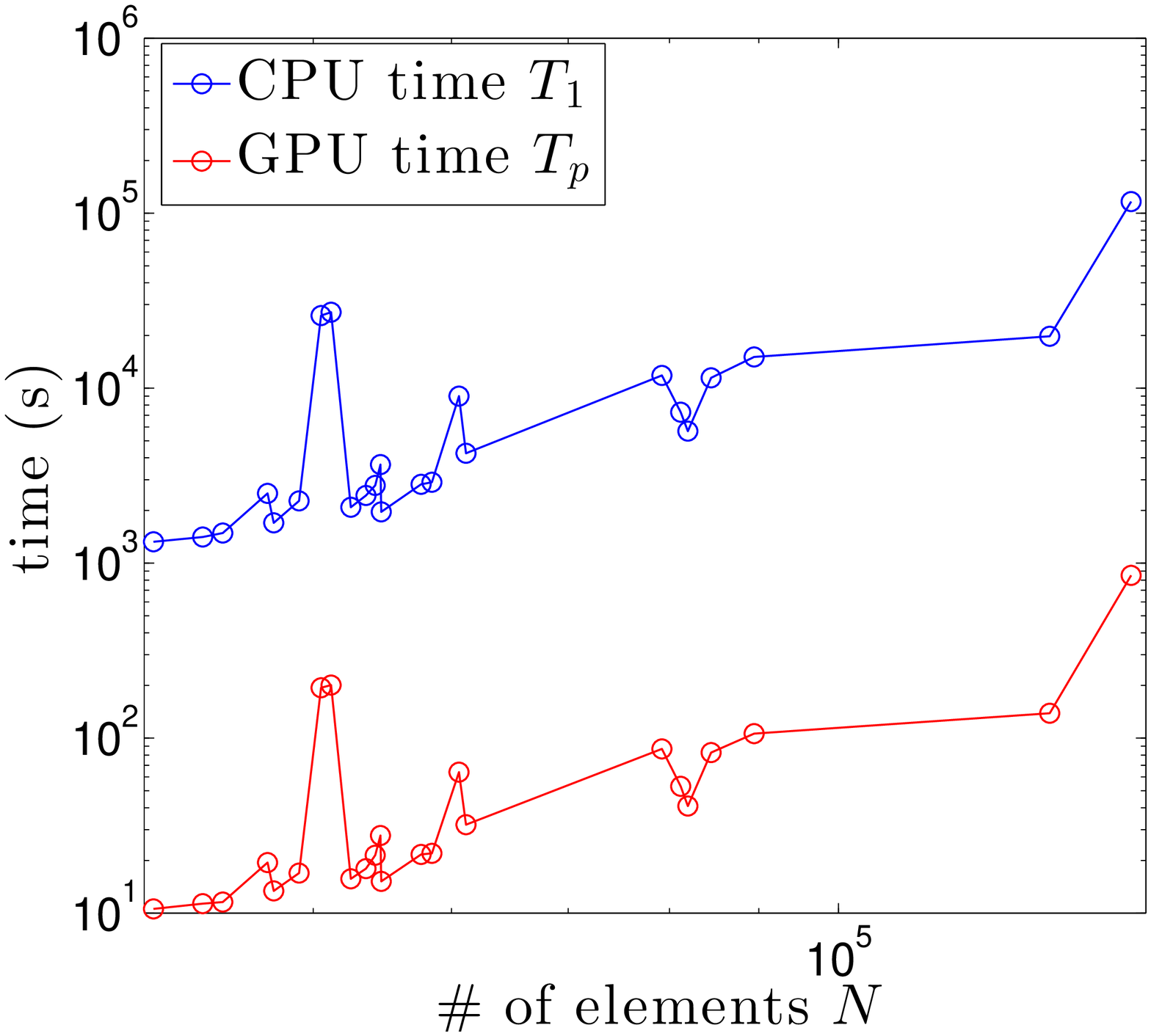}
\includegraphics[width=2.18in]{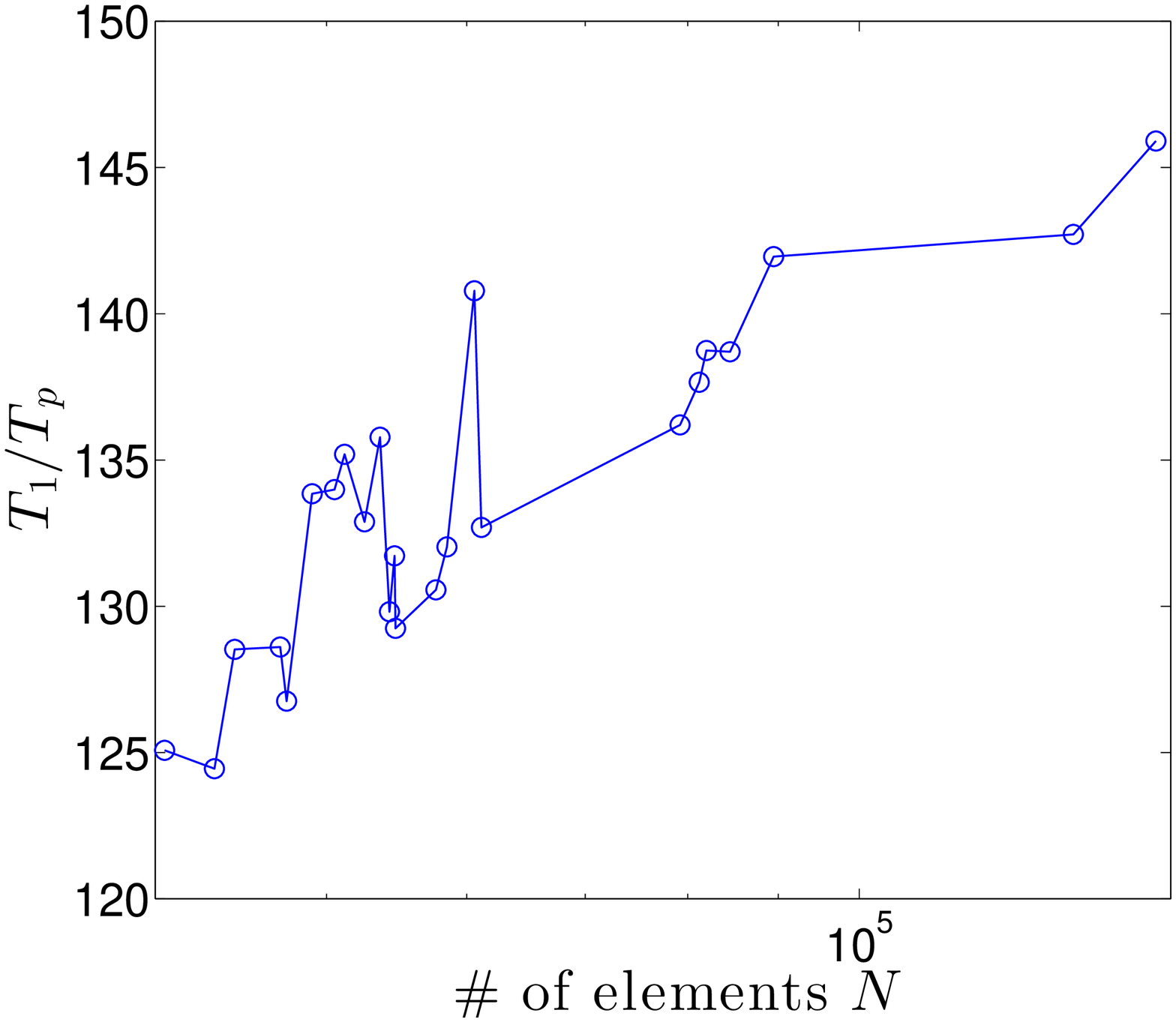}\\
\caption{Numerical results for solving PB equation and computing electrostatic solvation energy on a set of 24 proteins: (a) number of iterations; (b) CPU time $T_1$ and GPU time $T_p$; (c) the parallel speedup $T_1/T_p$.}
\label{fig_24p}
\end{center}
\end{figure}

Next we focus on two proteins: 1frd and 1svr. Results of these two proteins in Table~\ref{tb_24p} show fast convergence with number of iterations around 10 and they have surface areas more or less of 4500\AA$^2$ and number of elements is about 80,000 at ``density=10". To test the accuracy and efficiency of GABI-PB solver on these two proteins, we solve the PB equation and compute the electrostatic solvation energy at different ``densities"
ranging from 1 to 32 by doubling its value each time. The numerical results are shown in Table~\ref{tb_2p}.
We can see that as the densities are doubled each time, the number of elements are doubled approximately. The solvation energies at different densities are approaching to the values at the largest density 32. The running time $T_p$ with GPU acceleration  increases at the rate of $\mathcal{O}(N^2)$
after the number of elements are sufficiently large. Both cases show that, for solving PB equation and computing electrostatic solvation energy
on molecular surfaces discretized with nearly 300,000 elements, the time required is only about 10 minutes.
The number of iteration in Table \ref{tb_2p} shows that the finer resolution for the same molecular surface
will not increase the number of iterations.  To further investigate the convergence of accuracy on proteins, we plot the electrostatic solvation energy on Fig.~\ref{fig_2p}.  By using cubic interpolation, we can see the electrostatic solvation energy computed for both proteins eventually converge toward its interpolated value (the red ``$*$").

\begin{table}[htdp]
{\small
\caption{\small Numerical results for solving PB equation and computing the electrostatic solvation energy on protein 1frd and 1svr.}
\begin{center}
\begin{tabular}{c|cccc|cccc}
\hline\hline
&\multicolumn{4}{c|}{1frd}&\multicolumn{4}{c}{1svr}\\
\hline
density &\# of ele.&$E_\text{sol}$&GPU $T_p$(s) & \# of it. & \# of ele.&$E_\text{sol}$&GPU $T_p$(s) & \# of it. \\\hline
1	&10176	&-3360.69	&2	&18	&13974		&-2038.30	&2&11\\
2	&17710	&-3003.45	&2	&9	&22946		&-1833.41	&4&10\\
4	&34520	&-2931.05	&8	&9	&38116		&-1769.44	&10&9\\
8	&66294	&-2894.02	&28	&9	&71030		&-1735.37	&31&9\\
16	&134180	&-2880.90	&109	&9	&144168		&-1723.64	&140&10\\
32	&266702	&-2872.65	&423	&9	&284478		&-1716.94	&642&11\\
\hline
\end{tabular}
\end{center}
\label{tb_2p}
}
\end{table}%

\begin{figure}[!htbp]
\begin{center}
\flushleft \hskip 0.2in (a) \hskip 3.1in (b)\\
\includegraphics[width=3.3in]{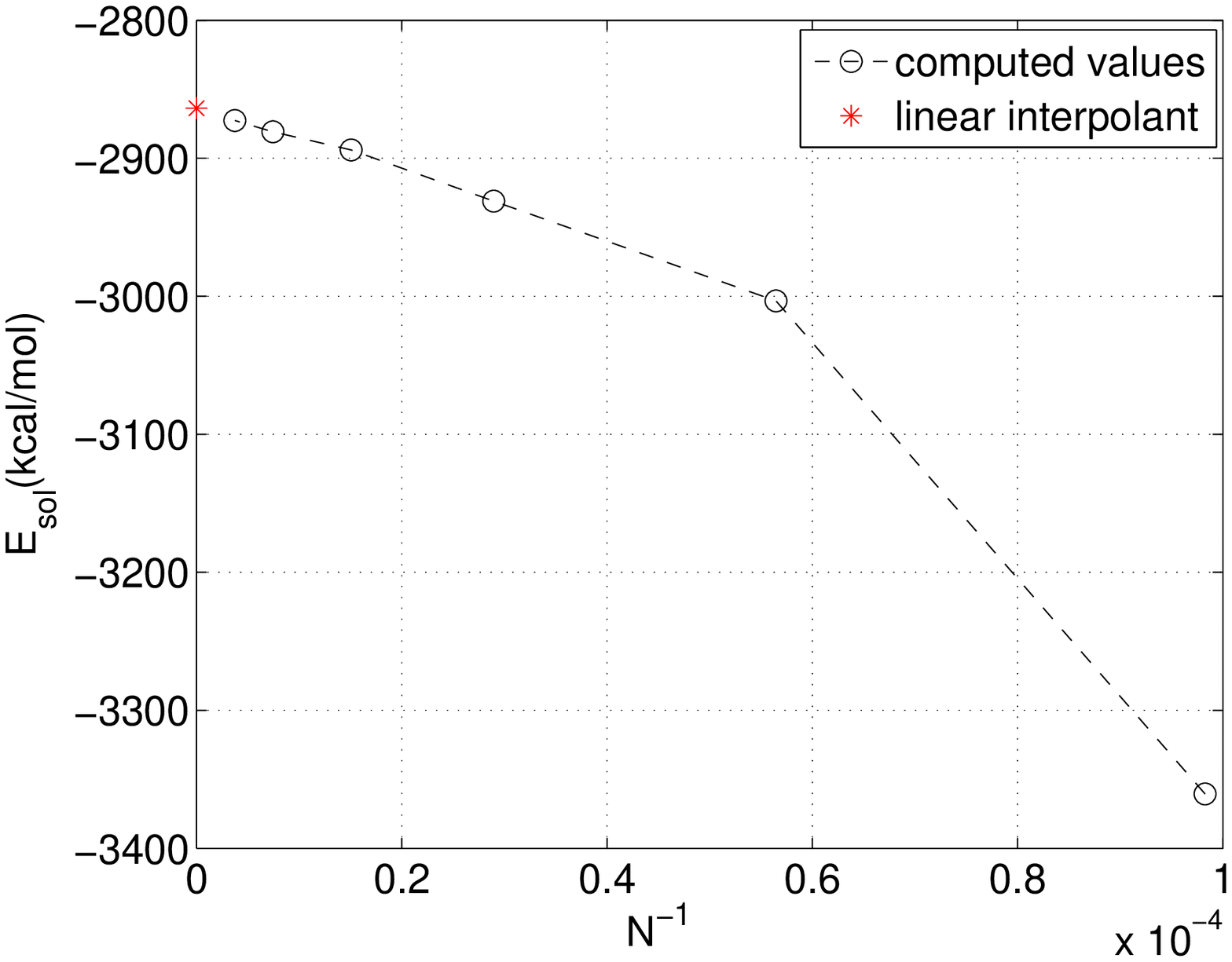}
\includegraphics[width=3.3in]{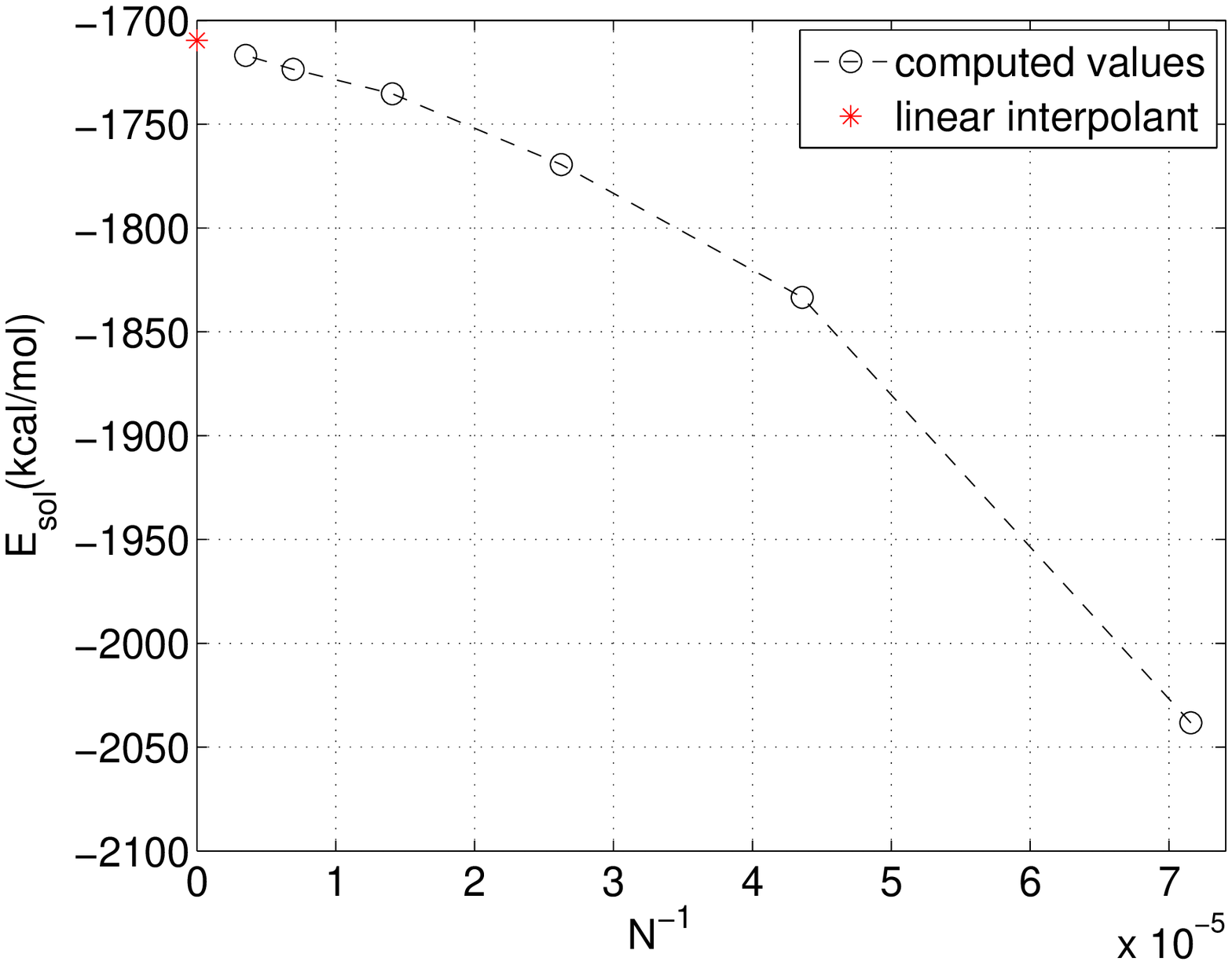}
\caption{Accuracy tests in terms for electrostatic solvation energy for proteins 1frd (a) and 1svr (b). }
\label{fig_2p}
\end{center}
\end{figure}

\section{Conclusion}
This paper describes a direct summation based GPU-accelerated boundary integral
Poisson-Boltzmann (GABI-PB) solver.
This solver discretizes the molecular surfaces with flat triangles and
performs numerical integration with centroid collocation schemes.
The numerical tests on spherical cavities show that GABI-PB solver
can achieve 0.5th order convergence on surface potentials
relative to the number of elements,
which is a 1st order convergence relative to the length.
The accurate surface potentials are of vital importance
to molecular modelings that are sensitive to electrostatics near or on the molecular surface. Meanwhile, 
this direct-sum boundary integral implementation uses memory efficiently. It has been shown that
we only need to allocate altogether $6N_c+11N$ size-of-double device memory,
where $N_c$ is the number of atoms and $N$ is the number of triangular elements.
In addition, numerical tests on a series of 24 proteins show fast convergence, consistent  electrostatic solvation energy computation, as well as 120-150X speedup using a single GPU (Nvidia Tesla M2070) card. 

The major limitation of the direct-sum scheme is obviously the $\mathcal{O}(N^2)$ computational cost, which becomes prohibitively expensive when the dimension of the problem increases to certain level.
Even though using multiple GPU cards can offset some of this effect, it costs extra design, implementation, and hardware purchase.
A remedy to this is applying the fast algorithms such as the $\mathcal{O}(N\log{N})$ treecode \cite{GengKrasnyjcp12} and the much more complicated and memory-consuming $\mathcal{O}(N)$ FMM \cite{Lu07}, which are under our investigation. It is hard to deny that the adoption of fast algorithms is necessary for large sized problems. However, there is a region on which the direct-sum GABI-PB solver has advantages over boundary integral PB solvers with fast algorithms. This region is bounded by what we called critical point where the fast algorithms surpasses the direct sum. For example, the CPU implementation of a treecode algorithm in solving boundary integral PB equation  \cite{GengKrasnyjcp12} shows a critical point at about $N=4000$ with $p=3$ (the order of Taylor expansion) and MAC $\leq 0.5$ (multipole acceptance criterion, the ratio between cluster radius and the distance between target particle and the center of the cluster). In GPU implementation, the critical point will be much bigger in considering the communication and memory access. We will identify the value of critical point in our future work. 

Memory usage is a critical factor of GPU performance. For solving boundary integral PB, direct-sum use about 1/3 of memory of treecode \cite{GengKrasnyjcp12} and 1/7 of memory of FMM \cite{Lu07}. The current GPU implementation is run on an available Tesla M2070 card. When we are considering to rerun all the tests on a Tesla M2090, Nvidia announced its release of Tesla K10 followed by K20 and K20x, showing the rapid hardware update year after year. Taking M2070 and K20x as examples for comparison, the peak double precision floating points performance increases from 515Gflops to 1.31Tflops, the memory bandwidth increases from 150 GB/sec to 250/sec, the number of CUDA cores increases from 448 to 2688. However the memory is still limited at 6Gbyte. These comparisons indicate the proposed direct-sum boundary integral PB solver, benefited from its low memory use, will demonstrate continuously improving performance with the flow of GPU hardware updates.  

In addition to including fast algorithms,
there are many spaces in which GABI-PB can be improved and extended.
For example, we are looking for better
triangulation programs for the molecular surfaces \cite{Bates08, Lujctc11, xu-zhang-09}
to avoid the slow convergence pattern for some proteins as seen in our tests.
A more challenging problem is the application of GABI-PB to molecular dynamics \cite{Gengjcp11, Lu03}, where the PB equation will be solved at every time sampling.
For molecular surfaces discretized within 50,000 elements, GABI-PB
can resolve each sample in a few seconds or less.
Furthermore, the
GPU-accelerated boundary integral scheme
has the potential to solve other integral equations such as the Helmholtz equation and Maxwell Equations.

\section*{Acknowledgements}
The work was supported by NSF grant DMS-0915057, University of Alabama new faculty startup fund
and Alabama Supercomputer Center. The authors thank Robert Krasny for helpful discussions.



\begin{thebibliography}{99}
\newcommand{\ncAddPaper} [7]{\bibitem{#1}#2,  {\it #3}, #4,    {\bf #5}, #6 (#7).}
\newcommand{\ncAddBook}  [5]{\bibitem{#1}#2,  {\it #3},  #4,    (#5).}
\newcommand{\ncAddPaperC}[6]{\bibitem{#1}#2, {\it  #3}, #4,    { #5}, (#6). }
\newcommand{\ncAddProced}[6]{\bibitem{#1}#2,  {\it #3},  in #4, {\bf #5}, #6.}



\ncAddPaper{BakerCurrOp} {N.A. Baker} {Improving implicit solvent
simulations: a Poisson-centric view} {Current Opin. Struct. Biol.}
{15} {137-143} {2005}




\ncAddPaper{Richards}{F.M. Richards}
   {Areas, volumes, packing and protein structure}
   {Annu. Rev. Biophys. Bioeng} {6} {151-176} {1977}

\ncAddPaper{Connolly85}{ M.L. Connolly} {Molecular surface
triangulation}    { J. Appl. Crystallogr.} {18} {499-505}{1985}

\ncAddPaper{Kirkwood34}
{J.G. Kirkwood}
{Theory of solution of molecules containing widely separated charges with special application to Zwitterions}
{J. Chem. Phys.} {7} {351-361} {1934}

\ncAddPaper{Holst} {M. Holst and F. Saied} {Multigrid solution of
the Poisson-Boltzmann equation} {J. Comput. Chem.} {14} {105-113}
{1993}

\ncAddPaper{Rocchia01} {W. Rocchia, E. Alexov and B. Honig}
{Extending the applicability of the nonlinear Poisson-Boltzmann
equation: multiple dielectric constants and multivalent ions} {J.
Phys. Chem. B} {105} {6507-6514} {2001}

\ncAddPaper{Im98}{W. Im, D. Beglov and B. Roux}
{Continuum solvation model: computtation of electrostatic forces
from numberical solutions to the Poisson-Boltzmann equation}
{Comp. Phys. Commun.}{111}{59-75}{1998}


\ncAddPaper{Luo02} {R. Luo, L. David and M.K. Gilson} {Accelerated
Poisson-Boltzmann calculations for static and dynamic systems} {J.
Comput. Chem.} {23} {1244-1253} {2002}


\ncAddPaper{bakerSept} {N.A. Baker, D. Sept, M.J. Holst and J.A.
McCammon} {The adaptive multilevel finite element solution of the
Poisson-Boltzmann
 equation on massively parallel computers}
{IBM J. Res. Dev.} {45} {427-438} {2001}

\ncAddPaper{Qiao}{Z. H. Qiao, Z.L. Li and T. Tang}{A finite
difference
 scheme for solving the nonlinear Poisson-Boltzmann equation
  modeling charged spheres}{J. Comp. Math.}{24}{252-264}{2006}

\ncAddPaper{Chenjcc11} {D. Chen, Z. Chen, C.J. Chen, W.H. Geng and G.W. Wei} {MIBPB: A Software Package for Electrostatic Analysis} {J. Comput. Chem.} {32} {756-770} {2011}


\ncAddPaper{JBKPB} {A. Juffer, E. Botta, B. van Keulen, A. van der
Ploeg and H. Berendsen} {The electric potential of a
macromolecule in a solvent: a fundamental approach} {J. Comput.
Phys.} {97} {144--171} {1991}

\ncAddPaper{BFZ} {A. Boschitsch, M.  Fenley and H.-X. Zhou} {Fast
boundary element method for the linear Poisson--Boltzmann
equation} {J. Phys. Chem.} {B 106} {2741-2754} {2002}

\ncAddPaper{Lu06}{B.Z. Lu, X.L. Cheng, J.F. Huang and J.A.  McCammon}
{Order N algorithm for computation of electrostatic interactions in biomolecular systems}
{PNAS} {103} {19314-19319} {2006}

\ncAddPaper{altman-bardhan-white-tidor-09}
{M.D. Altman, J.P. Bardhan, J.K. White, B. Tidor}
{Accurate solution of multi-region continuum biomolecule
electrostatic problems using the linearized Poisson-Boltzmann equation
with curved boundary elements}
{J. Comput. Chem. }{30}{132--153}{2009}


\ncAddPaper{yokota-bardhan-knepley-barba-hamada-11}
{R. Yokota, J.P. Bardhan, M.G. Knepley, L.A. Barba, T. Hamada}
{Biomolecular electrostatics using a fast multipole BEM on up to 512 GPUS and a billion unknowns}
{Comput. Phys. Commun.}{182} {1272--1283}{2011}

\ncAddPaper{GengKrasnyjcp12}
{W.H. Geng and R. Krasny}
{A Treecode-Accelerated Boundary Integral Poisson-Boltzmann Solver for
Electrostatics of Solvated Biomolecules}
{J. Comput. Phys.}  {submitted} {} {2012}

\ncAddPaper{Bordner} {A. Bordner and G. Huber} {Boundary element
solution of the linear Poisson--Boltzmann equation and a multipole
method for the rapid calculation of forces on macromolecules in
solution} {J. Comput. Chem.} {24} {353-367} {2003}

\ncAddPaper{LS} {J. Liang and S. Subranmaniam} {Computation of
molecular electrostatics with boundary element methods} {Biophys. J.} {73} {1830--1841} {1997}

\ncAddPaper{VorSch} {Y.N. Vorobjev and H.A. Scheraga} {A fast
adaptive multigrid boundary element method for macromolecular
electrostatic computations in a solvent} {J. Comput. Chem.} {18}
{569-583} {1997}

\ncAddPaper{bajaj-chen-rand-11}
{C. Bajaj, S.-C. Chen, A. Rand}
{An efficient higher-order fast multipole boundary element solution for Poisson-Boltzmann-based
molecular electrostatics}
{SIAM J. Sci. Comput.}{33}{826--848}{2011}

\ncAddPaper{AMBER_GPU_GB} {A. W. Goetz, M. J. Williamson, D. Xu, D. Poole, S. L. Grand and R. C. Walker}{Routine microsecond molecular dynamics simulations with AMBER - Part I: Generalized Born"}{}{in preparation}{}{2011}

\ncAddPaper{AMBER_GPU_PME}{A. W. Goetz, R. Salomon-Ferrer, D. Poole, S. L. Grand and R. C. Walker}{Routine microsecond molecular dynamics simulations with AMBER - Part II: Particle Mesh Ewald} {}{in preparation}{}{2011}

\ncAddPaper{NAMD}{J. C. Phillips, R. Braun, W. Wang, J. Gumbart, E. Tajkhorshid, E. Villa, C. Chipot, R. D. Skeel, L. Kale, and K. Schulten}{Scalable molecular dynamics with NAMD} {J. Comput. Chem}{26}{1781Ð1802}{2005}

\ncAddPaper{NAMD_GPU}
{J. E. Stone, D. J. Hardy, I. S. Ufimtsev, K. Schulten}{GPU-accelerated molecular modeling coming of age}{J. Molecular Graphics and Modelling}{29}{116-125}{2010}

\ncAddPaper{Nyland} {L. Nyland, M. Harris and J. Prins} {Fast N-Body Simulation
with CUDA} {GPU Gem 3} {Chapter 31} {677-695} {2009}

\ncAddPaper{Burtscher}{M. Burtscher and K. Pingali}{An Efficient CUDA Implementation of the Tree-based Barnes Hut n-Body Algorithm}
{GPU Computing Gems Emerald Edition} {Chapter 6} {75-92} {2011}

\ncAddPaper{barnes-hut-86} {J. Barnes and P. Hut}
{A hierarchical $O(N\log N)$ force-calculation algorithm}
{Nature} {324} {446-449} {1986}

\ncAddPaper{Yokota}{R. Yokota, J.P. Bardhan, M. G. Knepley, L.A. Barba, T. Hamada}
{Biomolecular electrostatics using a fast multipole BEM on up to 512 GPUs and a billion unknowns}
{Comput. Phys. Comm.}{182}{1272-1283}{2010}







%

%

\ncAddPaper{Yoon}
{B. J. Yong and A.M. Lenhoff}{A boundary element method for molecular electrostatics with electrolyte effects} {J. Comput. Chem.}{11}{1080-1086}{1990}

%


\ncAddPaper{Sanner} {M.F. Sanner, A.J. Olson and J.C. Spehner}
   {Reduced surface: An efficient way to compute molecular surfaces}
   {Biopolymers} {38} {305-320} {1996}


\ncAddPaper{Lu07}{B.Z. Lu, X.L. Cheng and J.A.  McCammon}
{''New-version-fast-multipole-method'' accelerated electrostatic calculations in biomolecular systems}
{J. Comput. Phys.} {226} {1348-1366} {2007}

\ncAddPaper{Saad} {Y. Saad and M.H. Schultz} {GMRES: A generalized minimal residual algorithm for solving nonsymmetric linear systems} {SIAM J. Sci. Stat. Comput} {7} {856-859} {1986}



%


\ncAddPaper{Gengjcp07} {W.H. Geng, S.N. Yu and G. W. Wei}
{Treatment of charge singularities in the implicit solvent
models} {J. Chem. Phys.}  {128} {114106} {2007}


\ncAddPaper{MacKerell-CHARMM22} {A.D. MacKerell Jr., D. Bashford,
M. Bellott, J.D. Dunbrack, M.J. Evanseck, M J. Field, S.
Fischer, J. Gao, H. Guo, S. Ha, D. Joseph-McCarthy, L. Kuchnir, K.
Kuczera, F.T.K. Lau, C. Mattos, S. Michnick, T. Ngo, D.T.
Nguyen, B. Prodhom, W.E. Reiher, B. Roux, M. Schlenkrich, J.C.
Smith, R. Stote, J. Straub, M. Watanabe, J. Wiorkiewicz-Kuczera,
D. Yin and M. Karplus} {All-atom empirical potential for molecular
modeling and dynamics studies of proteins} {J. Phys. Chem.} {102}
{3586-3616} {1998}


\ncAddPaper{Bates08} {P. Bates, G.W. Wei and S. Zhao} {Minimal
molecular  surfaces and their applications} {J.  Comput. Chem.}
{29} {380-391} {2008}

\ncAddPaper{Lujctc11}{M.X. Chen and B.Z. Lu}
{TMSmesh: A Robust Method for Molecular Surface Mesh
Generation Using a Trace Technique}
{J. Chem. Theory Comput.} {7} {203–212} {2011}

\bibitem{xu-zhang-09}
D. Xu and Y. Zhang,
{\it Generating triangulated macromolecular surfaces by Euclidean Distance Transform},
PLoS ONE 4(12): e8140. doi:10.1371/journal.pone.0008140

\ncAddPaper{GH}{L. Greengard and J. Huang}
{A new version of the fast multipole
method for screened Coulomb interactions in three dimensions}
{J. Comput. Phys.}{180}{642-658}{2002}

\ncAddPaper{Lijcp09}{P.J. Li, R. Krasny and H. Johnston}
{A Cartesian treecode for screened Coulomb particle interactions}
{J. Comput. Phys.} {228} {3858-3868} {2009}


\ncAddPaper{Gengjcp11}
{W.H. Geng and G.W. Wei}
{Multiscale molecular dynamics using the matched interface and boundary method}
{J. Comput. Phys.} {230} {435-457} {2011}

\ncAddPaper{Lu03} {Q. Lu and R. Luo} {A Poisson-Boltzmann dynamics
method with nonperiodic boundary condition} {J. Chem. Phys.} {119}
{11035-11047} {2003}



\end{thebibliography}
\end{document}